\documentclass[12pt,draft,leqno]{article}
\usepackage{amssymb, eucal, latexsym}

\newtheorem{theorem}{Theorem}[section]

\newtheorem{exa}[theorem]{Example}

\newtheorem{cor}[theorem]{Corollary}
\newtheorem{pro}[theorem]{Proposition}
\newtheorem{defi}[theorem]{Definition}

\newtheorem{rem}[theorem]{Remark}

\newtheorem{fact}[theorem]{Fact}

\def\a{\alpha}
\def\b{\beta}

\def\UP{\Upsilon}

\def\lra{\longrightarrow}

\def\sbe{\subseteq}

\def\stm{\setminus}
\def\ems{\emptyset}
\def\nes{\neq\emptyset}

\def\ap{^{\prime}}

\def\st{\ |\ }

\def\nin{\not\in}

\def\AA{{\cal A}}
\def\BB{{\cal B}}

\def\KK{{\cal K}}

\def\MM{{\cal M}}
\def\NN{{\cal N}}
\def\OO{{\cal O}}
\def\PP{{\cal P}}

\def\SS{{\cal S}}
\def\TT{{\cal T}}
\def\UU{{\cal U}}
\def\VV{{\cal V}}

\def\1{{\bf 1}}
\def\2{\mbox{{\bf 2}}}
\def\3{\mbox{{\bf 3}}}

\def\int{\mbox{{\rm int}}}

\def\doc{\hspace{-1cm}{\em Proof.}~~}
\def\sq{\hspace*{\fill} \hbox{\vrule\vbox{\hrule\phantom{o}\hrule}\vrule}}
\def\sqs{\sq \vspace{2mm}}

\def\RRRR{\mathbb{R}}
\def\NNNN{\mathbb{N}}

\def\tcx{t_X^C}
\def\tcy{t_Y^C}
\def\tcx0{t_{(X,X_0)}}
\def\tcy0{t_{(Y,Y_0)}}

\def\bU0{\bar{U}=(U^0,(U^i,U^{ci})_{i\in\omega})}

\def\bV0{\bar{V}=(V^0,(V^i,V^{ci})_{i\in\omega})}

\title{{\LARGE\bf
Lower-Vietoris-type Topologies on Hyperspaces}\\
\vspace{0.35cm}
{\large\bf Elza Ivanova-Dimova}\\
\vspace{0.25cm}
 {\footnotesize Dept. of Math. and
Informatics, Sofia University,  J. Bourchier 5, 1126 Sofia,
Bulgaria}}

\author{}

\date{}

\begin{document}
\maketitle
\begin{abstract}
{\footnotesize
We introduce a new
lower-Vietoris-type hypertopology in a way similar to that with which a new upper-Vietoris-type hypertopology was introduced in \cite{DV} (it was called there {\em Tychonoff-type hypertopology}). We study this new hypertopology and,
in particular, we generalize many results from \cite{CMR}.
As a corollary, we get that for every continuous map $f:X\longrightarrow X$, where $X$ is a continuum, there exist a subcontinuum $K$ of $X$ such that $f(K)=K.$
\noindent }
\end{abstract}

{\footnotesize {\em  MSC:} 54B20, 54H25, 54F15, 54C05.

}

\footnotetext[1]{{\footnotesize {\em E-mail address:}
elza@fmi.uni-sofia.bg}}

\baselineskip = \normalbaselineskip

\section{Introduction}

In 1975, M. M. Choban \cite{Ch} introduced a new topology on the set of all closed
subsets of a topological space for obtaining a generalization of the famous Kolmogoroff
Theorem on operations on sets. This new topology is similar to the upper Vietoris
topology but is weaker than it. In 1998, G. Dimov and D. Vakarelov \cite{DV} used a
generalized version of this new topology for proving an isomorphism theorem for the
category of all Tarski consequence systems. Later on it was studied in details in \cite{DOT}. In this paper we introduce a new
lower-Vietoris-type hypertopology in a way similar to that with which a new upper-Vietoris-type hypertopology was introduced in \cite{DV} (it was called there {\em Tychonoff-type hypertopology}). We study this new topology and,
in particular, we generalize many results from \cite{CMR}.

The paper is organized as follows. In Section 2 we introduce the notion of a {\it lower-Vietoris-type hypertopology}, we briefly study it and show that in general it differs from the lower Vietoris hypertopology (see Example
\ref{prm razl v-}). In this section, as well as in Sections 3, 4 and 5, we generalize many results from \cite{CMR}. In Section 3 we obtain a few results about some natural classes of maps between hyperspaces endowed with lower-Vietoris-type topologies. They generalize some results from \cite{CMR} and a theorem of H. J. Schmidt \cite[Theorem 11(1)]{Sch} about commutability between hyperspaces and subspaces. In Section 4, generalizing again some results from \cite{CMR}, we show when a hyperspace endowed with a lower-Vietoris-type topology is compact and when it has some given weight. In the last Section 5, generalizing some results from \cite{CMR}, we show that under mild conditions, the hyperspaces endowed with a lower-Vietoris-type topology have a trivial homotopy type, are absolute extensors for the class of all topological spaces and have the fixed-point property. As a corollary, we get that for every continuous map $f:X\longrightarrow X$, where $X$ is a continuum, there exist a subcontinuum $K$ of $X$ such that $f(K)=K.$

Let us fix the notation.

We denote by $\NNNN$ the set of all natural numbers (hence, $0\not\in\NNNN$), by $\RRRR$ the real line (with its natural topology) and by $\overline{\RRRR}$ the set $\RRRR\cup\{-\infty,\infty\}$.

Let $X$ be a set. We denote by $|X|$ the cardinality of $X$ and by $\PP(X)$ (resp., by $\PP\ap(X)$) the set of all (non-empty) subsets of $X$. Let $\MM,\AA\subseteq\PP(X)$ and $A\subseteq X$. We set:

$\bullet$ $A_{\MM}^+:=\{M\in\MM\st M\subseteq A\}$;

$\bullet$ $\AA_{\MM}^+:=\{A_{\MM}^+\st A\in\AA\}$;

$\bullet$ $A^-_{\MM}:=\{M\in\MM\st M\cap A\not=\emptyset\}$

$\bullet$ $\AA_{\MM}^-:=\{A_{\MM}^-\st A\in\AA\}$;

$\bullet$ ${\it Fin}(X):=\{M\subseteq X\st 0<|M|<\aleph_0\}$;

$\bullet$ ${\it Fin_n}(X):=\{M\subseteq X\st 0<|M|\leq n\}$, where $n\in\NNNN$;

$\bullet$ $\displaystyle\AA^{\cap}:=\{\bigcap_{i=1}^k A_i\st k\in\NNNN, A_i\in\AA\}$.


\noindent Let $(X,\TT)$ be a topological space. We put

$\bullet$ $CL(X):=\{M\subseteq X\st M \mbox{ is closed in }X,\ M\not=\emptyset\}$.

\noindent The closure of a subset $A$ of $X$ in $(X,\TT)$ will be denoted by $cl_XA$ or $\overline{A}^X$. When $\MM=CL(X)$, we will simply write $A^+$ and $A^-$ instead of $A_{\MM}^+$ and $A_{\MM}^-$.

All undefined here notions and notation can be found in \cite{AP,E2}.

\section{Lower-Vietoris-type topologies on hyperspaces}

Let $X$ be a topological space. Recall that the {\em  upper Vietoris topology} $\UP^+_X$ on $CL(X)$ (called also {\em Tychonoff topology on} $CL(X)$) has as a base the family of all sets of the form $$U_X^+=\{F\in CL(X)\st F\subseteq U\},$$ where $U$ is open in $X$, and {\em the lower Vietoris topology} $\UP_X^-$ on $CL(X)$ has as a subbase all sets of the form $$U_X^-=\{F\in CL(X)\st F\cap U\not=\emptyset\},$$ where $U$ is open in $X$.

\begin{defi}\label{Tichnoff type}{\rm (\cite{DV})}
\rm
(a) Let $(X,\TT)$ be a topological space and $\MM\subseteq\PP(X)$. The topology $\OO_{\TT}^+$ on $\MM$ having as a base the family $\TT_{\MM}^+$ will be called a {\em Tychonoff topology on} $\MM$ generated by $(X,\TT)$. When $\MM=CL(X)$, then $\OO_{\TT}^+$ is just the classical upper Vietoris topology on $CL(X)$.

\medskip

(b) Let $X$ be a set and $\MM\subseteq\PP(X)$. A topology $\OO$ on the set $\MM$ is called a {\em Tychonoff-type topology on} $\MM$ if the family $\OO\cap\PP(X)_{\MM}^+$ is a base for $\OO$.
\end{defi}

In what follows, the assertions whose proofs are (almost) obvious will be stated without any proofs.

\begin{fact}\label{fam-}
\rm
Let $X$ be a set, $\MM,\AA\subseteq\PP\ap(X)$. Then

a) $\bigcup\AA_{\MM}^-=(\bigcup\AA)_{\MM}^-$;

b) $A\subseteq B\Rightarrow A_{\MM}^-\subseteq B_{\MM}^-$.
\end{fact}

\begin{rem}\label{cap-}
\rm
Let us note that if $X$ is a set, $\MM\subseteq\PP\ap(X)$, ${\it Fin}_2(X)\subseteq\MM$, $A,B\subseteq X$ and $A\setminus B\not=\emptyset$, $B\setminus A\not=\emptyset$, then $(A\cap B)_{\MM}^-$ is a proper subset of $A_{\MM}^-\cap B_{\MM}^-$.

Indeed, it is obvious that $(A\cap B)_{\MM}^-\subseteq A_{\MM}^-\cap B_{\MM}^-$.
Let $x\in A\setminus B$ and $y\in B\setminus A$. Then $\{x,y\}\in\MM$ and $\{x,y\}\in A_{\MM}^-\cap B_{\MM}^-$, but $\{x,y\}\not\in (A\cap B)_{\MM}^-$.
\end{rem}

\begin{defi}\label{LVT}
\rm
Let $(X,\TT)$ be a topological space and $\MM\subseteq\PP\ap(X)$. The topology $\OO_{\TT}$ on $\MM$ having as a subbase the family $\TT_{\MM}^-$ will be called a {\em lower Vietoris topology on} $\MM$ generated by $(X,\TT)$. When $\MM=CL(X)$, then $\OO_{\TT}$ is just the classical lower Vietoris topology on $CL(X)$.
\end{defi}

\begin{defi}\label{lower-Vietoris-type}
\rm
Let $X$ be a set, $\MM\subseteq\PP\ap(X)$, $\OO$ be a topology on $\MM$. We say that $\OO$ is a {\em lower-Vietoris-type topology on} $\MM$, if $\OO\cap\{A^-_{\MM}\st A\subseteq X\}$ is a subbase for $\OO$.
\end{defi}

Clearly, a lower Vietoris topology on $\MM$ is always a lower-Vietoris-type topology on $\MM$, but not viceversa (see Example \ref{prm razl v-}).

\begin{fact}\label{Vietoris-type}
\rm
Let $X$ be a set, $\MM\subseteq\PP\ap(X)$ and $\OO$ be a lower-Vietoris-type topology on $\MM$. Then the family $$\PP_{\OO}:=\{A\subseteq X\st A_{\MM}^-\in\OO\}$$ contains $X$, and, hence, can serve as a subbase of a topology $$\TT_{\OO}$$ on $X$. The family $(\PP_{\OO})_{\MM}^-$ is a subbase for $\OO$. The family $\PP_{\OO}$ is closed under arbitrary unions.
\end{fact}

\begin{defi}\label{2.5}
\rm
Let $X$ be a set, $\MM\subseteq\PP\ap(X)$ and $\OO$ be a lower-Vietoris-type topology on $\MM$. Then we will say that the topology $\TT_{\OO}$ on $X$, introduced in Fact \ref{Vietoris-type}, is {\em induced by the topological space} $(\MM,\OO)$.
\end{defi}

\begin{pro}\label{subbase1}
Let $X$ be a set and $\MM\subseteq\PP\ap(X)$. Then a topology $\OO$ on $\MM$ is a lower-Vietoris-type topology on $\MM$ iff there exists a topology $\TT$ on $X$ and a subbase $\SS$ for $\TT$ (which contains  $X$ and is closed under arbitrary unions), such that $\SS^-_{\MM}=\{A_{\MM}^-\st A\in\SS\}$ is a subbase for $\OO$.
\end{pro}

\begin{pro}\label{subbase2}
Let $X$ be a set, $\MM\subseteq\PP\ap(X)$ and $\SS\subseteq\PP\ap(X)$. Then $\SS^-_{\MM}$ is a subbase for a (lower-Vietoris-type) topology on $\MM$
if and only if $M\cap\bigcup\SS\not=\emptyset$ for any $M\in\MM$.
\end{pro}

\begin{defi}\label{gene}
\rm
Let $X$ be a set, $\MM,\SS\subseteq\PP\ap(X)$ and $M\cap\bigcup\SS\not=\emptyset$ for any $M\in\MM$. Then the lower-Vietoris-type topology $$\OO_\SS^\MM$$ on $\MM$
for which $\SS^-_\MM$ is a subbase (see Proposition \ref{subbase2}) will be called {\em lower-Vietoris-type topology on $\MM$ generated by the family $\SS$.}
When there is no ambiguity, we will simply write $$\OO_\SS$$ instead of $\OO_\SS^\MM$.
\end{defi}

\begin{cor}\label{subbase3}
Let $X$ be a set, $\MM,\SS\subseteq\PP\ap(X)$ and $\bigcup\SS=X$. Then $\SS^-_{\MM}$ is a subbase for a lower-Vietoris-type topology $\OO_\SS$ on $\MM$.
\end{cor}

\begin{pro}\label{equiv-top}
Let $X$ be a set, $\MM\subseteq\PP\ap(X)$ and $\UU_i\subseteq\PP\ap(X)$, $i=1,2$, be such that $X=\bigcup\UU_1=\bigcup\UU_2$. Let $\OO_i$ be the topology on $\MM$ generated by $\UU_i$, $i=1,2$. Then $\OO_1\equiv\OO_2$ if and only if $\PP_{\OO_1}\equiv\PP_{\OO_2}$.
\end{pro}

\doc Clearly, if $\OO_1\equiv\OO_2$ then $\PP_{\OO_1}\equiv\PP_{\OO_2}$. Conversely, let $\PP_{\OO_1}\equiv\PP_{\OO_2}$. Obviously, $\UU_i\subseteq\PP_{\OO_i}$ for $i=1,2$. Thus $(\PP_{\OO_i})^-_{\MM}$ is a subbase of $\OO_i$, $i=1,2$. Hence $\OO_1\equiv\OO_2$. \sqs

\begin{fact}\label{zatwob}
Let $X$ be a set, $\MM\subseteq\PP\ap(X)$  and $\OO$ be a lower-Vietoris-type topology on $\MM$. Then, for every $C\in\MM$,
we have that $\overline{C}^{(\MM,\OO)}\supseteq\{M\in\MM\st M\subseteq C\}=C^+_\MM$ .
\end{fact}

\begin{pro}\label{V^-}
If $(X,\TT)$ is a $T_1$-space, $\MM=CL(X)$ (or $\MM$ is a closed base for $X$), $\OO$ is a lower-Vietoris-type topology on $\MM$ and $\TT=\TT_{\OO}$ (see Fact \ref{Vietoris-type} for $\TT_\OO$), then $\OO\equiv \UP^-_X$ if and only if for every $F\in\MM$ we have that $\overline{F}^{\OO}=F^+$.
\end{pro}

\doc
It is obvious that $\overline{F}^{\UP^-_X}=F^+$, for every $F\in\MM$.

Conversely, let $\overline{F}^{\OO}=F^+$ for every $F\in\MM$. By Fact \ref{Vietoris-type}, for proving that $\OO\equiv\UP^-_X$,
it suffices to show that $\PP_{\OO}$ is a base for $\TT$.
Let $F$ be a closed subset of $X$ and $x\not\in F$. Then $\{x\}\not\in\overline{F}^{\OO}$. Hence there exist $U_1,\dots,U_n\in\PP_{\OO}$, such that $\displaystyle\{x\}\in\bigcap_{i=1}^nU_i^-=U$ and $F\not\in U$. Then there exists an $i\in\{1,\dots,n\}$ such that $F\not\in U_i^-$, i.e. $F\cap U_i=\emptyset$. Thus $x\in U_i\subseteq X\setminus F$ and $U_i\in\PP_\OO$. Hence $\PP_\OO$ is a base for $\TT$. Therefore $\PP_{\OO}\equiv\TT$, i.e. $\OO\equiv \UP^-_X$ .\sqs

\begin{cor}\label{kuratovski}
 Let $(X,\TT)$ be a $T_1$-space, $\MM=CL(X)$ (or $\MM$ is a closed base for $X$), $\OO$ is a lower-Vietoris-type topology on $\MM$ and $\TT=\TT_{\OO}$. Then the Kuratowski operators for $\OO$ and $\UP^-_X$ coincide on the singletons of $\MM$ iff  they coincide on every subset of $\MM$.
\end{cor}

\begin{pro}\label{t0}
Let $X$ be a set, $\MM\subseteq\PP\ap(X)$, $\OO$ be a lower-Vietoris-type topology on $\MM$. Then $(\MM,\OO)$ is a $T_0$-space if and only if for every $M_1,M_2\in\MM$ with $M_1\not=M_2$, there exist $U_1,\dots,U_n\in\PP_{\OO}$ such that either ($M_1\cap U_i\not=\emptyset$ for any $i=1,\dots,n$ and there exists an $i_0\in\{1,\dots,n\}$ for which $M_2\cap U_{i_0}=\emptyset$) or ($M_2\cap U_i\not=\emptyset$ for any $i=1,\dots,n$ and there exists an $i_0\in\{1,\dots,n\}$ for which $M_1\cap U_{i_0}=\emptyset$).
\end{pro}

\begin{defi}\label{natural}
\rm
Let $X$ be a set and $\MM\subseteq\PP(X)$. We say that $\MM$ is a {\em natural family in} $X$ if $\{x\}\in\MM$ for any $x\in X$.
\end{defi}

\begin{pro}\label{t1}
Let $X$ be a set, $\MM\subseteq\PP\ap(X)$ be a natural family, $\OO$ be a lower-Vietoris-type topology on $\MM$. Then $(\MM,\OO)$ is a $T_1$-space if and only if $$\MM=\{\{x\}\st x\in X\}$$ and $(X,\TT_{\OO})$ is an $T_1$-space.
\end{pro}

We are now going to construct a lower-Vietoris-type topology on a hyperspace, which is not a lower Vietoris topology. For doing this we will need some preliminary definitions and statements.

\begin{pro}\label{U-}
Let $X$ be a set, $\MM\subseteq\PP\ap(X)$ be a natural family,  $$\{U_{\a,i}\st \a\in A,i=1,\ldots,n_\a, n_\a\in\NNNN\}\sbe\PP\ap(X)$$ and $U\sbe X$. Then $\displaystyle U^-_\MM=\bigcup_{\alpha\in A}\bigcap_{i=1}^{n_{\alpha}}(U_{\alpha,i})^-_\MM$ if and only if the next two conditions are fulfilled:

\smallskip

\noindent(1) $\displaystyle U=\bigcup_{\alpha\in A}\bigcap_{i=1}^{n_{\alpha}}U_{\alpha,i}$, and

\smallskip

\noindent(2) for every $M\in\MM$ such that $M\cap U_{\alpha,i}\not=\emptyset$ for some $\alpha\in A$ and for any $i=1,\dots,n_{\alpha}$, we have that $M\cap U\not=\emptyset$.
\end{pro}

\doc Let $\displaystyle U^-_\MM=\bigcup_{\alpha\in A}\bigcap_{i=1}^{n_{\alpha}}(U_{\alpha,i})^-_\MM$.
 For every $x\in X$,
we have that ($x\in U$) $\iff$ ($\{x\}\in U^-_{\MM}$) $\iff$ (there exists an $\alpha\in A$ with $\displaystyle x\in\bigcap_{i=1}^{n_{\alpha}}U_{\alpha,i}$) $\iff$ ($\displaystyle x\in\bigcup_{\alpha\in A}\bigcap_{i=1}^{n_{\alpha}}U_{\alpha,i}$). Hence $\displaystyle U=\bigcup_{\alpha\in A}\bigcap_{i=1}^{n_{\alpha}}U_{\alpha,i}\in\TT_{\OO}$. So, condition (1) is fulfilled. Clearly, condition (2) is also fulfilled.

Conversely, let $M\in U^-_\MM$. Then $M\in\MM$ and there exists  $x\in M\cap U$. Hence by condition (1), there exists  $\alpha\in A$ such that $\displaystyle x\in\bigcap_{i=1}^{n_{\alpha}}U_{\alpha,i}$. Thus $\displaystyle M\in\bigcup_{\alpha\in A}\bigcap_{i=1}^{n_{\alpha}}(U_{\alpha,i})^-_\MM$. Therefore $\displaystyle U^-_\MM=(\bigcup_{\alpha\in A}\bigcap_{i=1}^{n_{\alpha}}U_{\alpha,i})^-_\MM\subseteq\bigcup_{\alpha\in A}\bigcap_{i=1}^{n_{\alpha}}(U_{\alpha,i})^-_\MM$.
Now, let $\displaystyle M\in\bigcup_{\alpha\in A}\bigcap_{i=1}^{n_{\alpha}}(U_{\alpha,i})^-_\MM$. Then there exists  $\alpha\in A$ such that $M\cap U_{\alpha,i}\not=\emptyset$ for any $i=1,\dots,n_{\alpha}$. Hence by condition (2), $M\cap U\nes$ and thus $M\in U^-_\MM$. \sqs

\begin{defi}\label{m^- cl}
\rm
Let $X$ be a set, $\MM\subseteq\PP\ap(X)$ be a natural family in $X$ and $\UU\subseteq\PP\ap(X)$. Then:

\smallskip

\noindent(a) If $\UU\ap=\{U_{\a,i}\st \a\in A,i=1,\ldots,n_\a, n_\a\in\NNNN\}$ is a subfamily of $\UU$,  $\displaystyle U=\bigcup_{\alpha\in A}\bigcap_{i=1}^{n_{\alpha}}U_{\alpha,i}$ and  from $M\in\MM$ and $M\cap U_{\alpha,i}\not=\emptyset$ for some $\alpha\in A$ and for any $i=1,\dots,n_{\alpha}$, it follows that $M\cap U\not=\emptyset$, then we will say that {\em the set $U$ is $\MM^-$-covered by the family $\UU\ap$;}

\smallskip

\noindent(b)
 The family $\UU$ is said to be an $\MM^-$-{\em closed family} if it contains any subset $U$ of $X$ which is $\MM^-$-covered by some subfamily $\UU\ap$ of $\UU$.
\end{defi}

\begin{pro}\label{m-cl}
Let $X$ be a set, $\MM\subseteq\PP\ap(X)$, $A$ be an index set, $\UU_{\alpha}\subseteq\PP\ap(X)$ and $\UU_{\alpha}$ are $\MM^-$-closed families for any $\alpha\in A$. Then $\UU=\bigcap\{\UU_{\alpha}\st \alpha\in A\}$ is an $\MM^-$-closed family.
\end{pro}

\doc Let $\displaystyle U=\bigcup_{\beta\in B}\bigcap_{i=1}^{n_{\beta}}U_{\beta,i}$, where $U_{\beta,i}\in\UU$ for any $\beta\in B$ and any $i=1,\ldots,n_\b$, and for every $M\in\MM$, such that there exists  $\beta\in B$ with $M\cap U_{\beta,i}\not=\emptyset$ for any $i=1,\dots,n_{\beta}$, we have that $M\cap U\not=\emptyset$. Then $U\in\UU_{\alpha}$ for any $\alpha\in A$. Hence $U\in\UU$. \sqs

\begin{cor}\label{m-zatwcor}
Let $X$ be a set and $\MM\subseteq\PP\ap(X)$. Then every
 family $\NN\subseteq\PP\ap(X)$ is contained in a minimal $\MM^-$-closed family, denoted by $$\MM^-(\NN).$$
\end{cor}

\doc Clearly, $\PP'(X)$ is an $\MM^-$-closed family. Now we can apply Proposition \ref{m-cl}. \sqs

\begin{pro}\label{P_O}
Let $X$ be a set, $\MM\subseteq\PP\ap(X)$ be a natural family and $\OO$ be a lower-Vietoris-type topology on $\MM$.  Then $\PP_{\OO}$ (where $\PP_{\OO}=\{A\subseteq X\st A^-_\MM\in\OO\}$)  is an $\MM^-$-closed family. If $\UU\sbe\PP\ap(X)$, $\UU$ covers $X$ and generates $\OO$, $\VV\sbe\PP\ap(X)$ and $\UU\sbe\VV\sbe\PP_\OO$, then
$\PP_\OO=\MM^-(\VV)$.
\end{pro}

\doc Let $\displaystyle U=\bigcup_{\alpha\in A}\bigcap_{i=1}^{n_{\alpha}}U_{\alpha,i}$, where $U_{\alpha,i}\in\PP_{\OO}$, $\forall\alpha\in A$ and $\forall i\in\{1,\dots,n_{\alpha}\}$ ($n_{\alpha}\in\NNNN$, $\forall\alpha\in A$). Let, also, for every $M\in\MM$ such that $M\cap U_{\alpha,i}\not=\emptyset$ for some $\alpha\in A$ and every $i\in\{1,\dots, n_{\alpha}\}$, we have that $M\cap U\not=\emptyset$. Then, by Proposition \ref{U-}, $\displaystyle U_{\MM}^-=\bigcup_{\alpha\in A}\bigcap_{i=1}^{n_{\alpha}}(U_{\alpha,i})_{\MM}^-$. Hence $U_{\MM}^-\in\OO$. This implies that $U\in\PP_{\OO}$. So, $\PP_{\OO}$ is an $\MM^-$-closed family. If $\UU'\subseteq\PP'(X)$ is an $\MM^-$-closed family and $\UU'\supseteq\VV$, then arguing as above, we get that $\PP_{\OO}\subseteq\UU'$. Hence $\MM^-(\VV)=\PP_{\OO}$. \sqs

\begin{cor}\label{sl1}
Let $X$ be a set, $\MM\subseteq\PP\ap(X)$ be a natural family and $\UU,\VV\subseteq\PP\ap(X)$ are covers of $X$. Then $\OO_{\UU}^\MM\equiv\OO_{\VV}^\MM$ if and only if $\MM^-(\UU)=\MM^-(\VV)$.
\end{cor}

\doc It follows from Propositions \ref{equiv-top} and \ref{P_O}. \sqs

\begin{cor}\label{sl2}
Let $(X,\TT)$ be a topological space, $\MM\subseteq\PP\ap(X)$ be a natural family and $\BB$ be a base for $\TT$. Then $\OO_{\BB}^\MM\equiv \OO^\MM_\TT$.
\end{cor}

\doc Set $\OO=\OO_{\BB}^{\MM}$. By Fact \ref{Vietoris-type}, $\PP_{\OO}$ is closed under arbitrary unions. Since $\BB\subseteq\PP_{\OO}$, we get that $\BB\subseteq\TT\subseteq\PP_{\OO}$. Thus, by Proposition \ref{P_O}, $\MM^-(\BB)=\PP_{\OO}=\MM^-(\TT)$. Now Corollary \ref{sl1} implies that $\OO_{\BB}^{\MM}\equiv\OO_{\TT}^{\MM}$. \sqs

\begin{exa}\label{prm razl v-}
\rm
Let $\MM\subseteq\PP\ap(\RRRR)$ and $\MM\supseteq Fin_2(\RRRR)$. Then $$\UU=\{(-\infty,\beta),(\alpha,\infty),(-\infty,\beta)\cup(\alpha,\infty)\st\alpha,\beta\in\overline{\RRRR}\}$$ is an $\MM^-$-closed family. It generates a topology $\OO_\UU^\MM$ on $\MM$ different from the lower Vietoris topology $\OO^\MM_\TT$ (= $\OO_\TT$) on $\MM$, where $\TT$ is the natural topology on $\RRRR$.
\end{exa}

\doc Let $\displaystyle U=\bigcup_{\alpha\in A}\bigcap_{i=1}^{n_{\alpha}}U_{\alpha,i}$, where $U_{\alpha,i}\in\UU$ for any $\a\in A$ and any $i\in\{1,\ldots,n_\a\}$, and let for every $M\in\MM$ such that there exists an $\alpha\in A$ with $M\cap U_{\alpha,i}\not=\emptyset$ for any $i=1,\dots,n_{\alpha}$, we have that $M\cap U\not=\emptyset$. We will prove that $U\in\UU$. Since $U\in\TT$, we have that $\displaystyle U=\biguplus_{i=1}^{\infty}(x_i,y_i)$ ($x_i,y_i\in\overline{\RRRR}$, $x_i<y_i$). Suppose that $U\not\in\UU$. Then there exists an $i\in\NNNN$ with $x_i,y_i\in\RRRR$  and $x_i,y_i\not\in U$. Let $M=\{x_i,y_i\}$.
Then $M\in\MM$. We have that $\displaystyle(x_i,y_i)\subseteq\bigcup_{\alpha\in A}\bigcap_{j=1}^{n_{\alpha}}U_{\alpha,j}$. Let $c\in(x_i,y_i)$. Then there exists an $\alpha\in A$ with $\displaystyle c\in\bigcap_{j=1}^{n_{\alpha}}U_{\alpha,j}$. Let $j\in\{1,\dots,n_{\alpha}\}$. If $U_{\alpha,j}=(-\infty,u)$, then $u>c$. Hence $x_i\in U_{\alpha,j}$, i.e. $M\cap U_{\alpha,j}\not=\emptyset$. If $U_{\alpha,j}=(v,\infty)$ then $v<c$. Hence $y_i\in(v,\infty)$, i.e. $M\cap U_{\alpha,j}\not=\emptyset$. If $U_{\alpha,j}=(-\infty,u)\cap(v,\infty)$, then $c\in(-\infty,u)$ or $c\in(v,\infty)$ and we get as above  that $M\cap U_{\alpha,j}\not=\emptyset$. Thus $M\cap U_{\alpha,j}\not=\emptyset$ for any $j\in\{1,\dots,n_{\alpha}\}$, and  $M\cap U=\emptyset$. We get a contradiction. Hence $U\in\UU$.
Now all follows from Proposition \ref{P_O}. \sqs

The next two propositions generalize, respectively, Corollary 1.3(b) and Proposition 1.4 from \cite{CMR}.

\begin{pro}\label{sepber}
Let $X$ be a set, $\MM\subseteq\PP\ap(X)$, $X\in\MM$ and $\OO$ be a lower-Vietoris-type topology on $\MM$. Then $(\MM,\OO)$ is a pseudocompact, connected and separable space. Also, the intersection of any family of open dense subsets of $(\MM,\OO)$ is a dense subset of $(\MM,\OO)$.
\end{pro}

\doc We have that $\overline{X}^{\OO}=\MM$. Hence $(\MM,\OO)$ is an separable space and every continuous mapping $f:(\MM,\OO)\longrightarrow Y$, where $Y$ is a $T_1$-space, is a constant map (indeed, $f(\overline{X}^{(\MM,\OO)})\subseteq\overline{f(X)}$, i.e. $|f(\MM)|=1$). Thus $(\MM,\OO)$ is pseudocompact and connected.

Let for any $\alpha\in A$, $\UU_{\alpha}\subseteq\MM$ be an open and dense subset of $(\MM,\OO)$. We will prove that $\UU=\bigcap\{\UU_{\alpha}\st\alpha\in A\}$ is a dense subset of $(\MM,\OO)$. First, we will prove that for any $\alpha\in A$, $X\in\UU_{\alpha}$. Indeed, let $F\in\UU_{\alpha}$. Then there exist $V_1,\dots,V_n\in\PP_\OO$ with $\displaystyle F\in\bigcap_{j=1}^n(V_j)^-_\MM\subseteq\UU_{\alpha}$. Since, for any $j=1,\ldots,n$, $X\cap V_j\not=\emptyset$, we have that $\displaystyle X\in\bigcap_{j=1}^n(V_j)^-_\MM$. Hence $X\in\UU_{\alpha}$.
So, $X\in\bigcap\{\UU_{\alpha}\st\alpha\in A\}=\UU$. Since $\overline{X}^{\OO}=\MM$, we get that $\UU$ is dense in $(\MM,\OO)$. \sqs

The following fact is obvious.

\begin{fact}\label{neprpr}
Let $(X,\TT)$ and $(Y,\TT')$ be topological spaces. Then $f:(X,\TT)\longrightarrow(Y,\TT')$ is continuous if and only if there exists a subbase $\PP$ of $\TT'$ such that if $x\in X$, $U\in\PP$ and $f(x)\in U$, then there exists a $V\in\TT$ such that $x\in V$ and $f(V)\subseteq U$.
\end{fact}

\begin{pro}\label{vl}
Let $X$ be a set, $\MM\subseteq\PP\ap(X)$, $\MM$ be a natural family and $\OO$ be a lower-Vietoris-type topology on $\MM$. Then
$\Phi:(X,\TT_{\OO})\longrightarrow(\MM,\OO)$, where $\Phi(x)=\{x\}$ for any $x\in X$, is a homeomorphic embedding. If $(X,\TT_{\OO})$ is a $T_2$-space then $\Phi(X)$ is a closed subset of $(\MM,\OO)$ and if, in addition, $|X|>1$ and $X\in\MM$, then $\Phi(X)$ is nowhere dense in $(\MM,\OO)$.
\end{pro}

\doc Let $x\in X$, $U\in\PP_{\OO}$ and $\Phi(x)=\{x\}\in U^-_\MM$. Then for any $y\in U$ we have that $\Phi(y)=\{y\}\in U^-_\MM$. Hence $\Phi$ is continuous. Conversely, let $x\in U\in\TT_{\OO}$. Then there exist $U_1,\dots,U_n\in\PP_\OO$ with $\displaystyle x\in\bigcap_{i=1}^nU_i\subseteq U$. Let $\displaystyle\{y\}\in\bigcap_{i=1}^n(U_i)^-_\MM$. Then $\displaystyle y\in\bigcap_{i=1}^nU_i\subseteq U$. Hence $\Phi$ is a homeomorphic embedding.

We will prove that $\Phi(X)$ is closed in $(\MM,\OO)$. Indeed, let $M\in\MM\setminus\Phi(X)$. Then $|M|\geq 2$. Since $(X,\TT_{\OO})$ is a $T_2$-space, we get that there exist $U_i\in\PP_\OO$, $i=1,\dots,n$ and $V_j\in\PP_\OO$, $j=1,\dots,m$ such that if $\displaystyle U=\bigcap_{i=1}^nU_i$ and $\displaystyle V=\bigcap_{j=1}^mV_j$ then $U\cap V=\emptyset$, $U\cap M\not=\emptyset$ and $V\cap M\not=\emptyset$. Then $\displaystyle M\in U^-_\MM\cap V^-_\MM=O\in\OO$ and $O\cap\Phi(X)=\emptyset$ (indeed, if $\{x\}\in O$ then $x\in U\cap V=\emptyset$,  a contradiction).

If $|X|>1$, then  $X\in\MM\setminus\Phi(X)$ and $\overline{X}^{(\MM,\OO)}=\MM$. Thus we get that $\Phi(X)$ is nowhere dense in $(\MM,\OO)$. \sqs

\begin{pro}\label{fingyst}
Let $X$ be a set, $\MM\subseteq\PP\ap(X)$, $\OO$ be a lower-Vietoris-type topology on $\MM$
and $\MM\supseteq Fin(X)$. Then $Fin(X)$ is dense in $(\MM,\OO)$.
\end{pro}

\doc Let $U_1,\dots, U_n\in\PP_{\OO}$ and $\displaystyle U=\bigcap_{i=1}^n(U_i)_{\MM}^-\not=\emptyset$. Let $M\in U$. Then, $\forall i\in\{1,\dots,n\}$, there exists an $x_i\in M\cap U_i$. Let $F=\{x_1,\dots,x_n\}$. Then $F\in Fin(X)\cap U$. So, $Fin(X)$ is dense in $(\MM,\OO)$. \sqs

\section{Maps between hyperspaces endowed with a lo\-wer-Vietoris-type topology}

With the next three results we generalize, respectively, Propositions 2.1, 2.2 and 2.3 from \cite{CMR}.

\begin{pro}\label{2^f}
Let $(X,\TT)$ and $(X',\TT')$ be two topological spaces,  $f:(X,\TT)\longrightarrow(X',\TT')$, $\PP$ be a subbase for $\TT$, $\PP\ap$ be a subbase for $\TT'$, $f^{-1}(\PP')\subseteq\PP$, $\MM\subseteq\PP'(X)$, $\MM'=CL(X')$. Let $\OO=\OO_{\PP}^{\MM}$ and $\OO'=\OO_{\PP'}^{\MM'}$. Then the map $2^f:(\MM,\OO)\longrightarrow(\MM',\OO')$, where $2^f(C)=\overline{f(C)}^{X'}$, for every $C\in \MM$, is continuous.
\end{pro}

\doc Let $U'\in\PP'$. Then

\noindent$\begin{array}{c}
(2^f)^{-1}((U')^-_{\MM\ap})=\{F\in\MM\st\overline{f(F)}^{X\ap}\cap U'\not=\emptyset\}=\{F\in\MM\st f(F)\cap U'\not=\emptyset\}= \\
\{F\in\MM\st F\cap f^{-1}(U')\not=\emptyset\}=(f^{-1}(U'))^-_\MM\in\OO.
\end{array}$

Hence $2^f$ is continuous. \sqs

\begin{pro}\label{komp}

\noindent(a) Let $(X,\TT)$ be a topological space, $\MM=CL(X)$, $\PP$ be a subbase for $\TT$ and $\OO=\OO_{\PP}^{\MM}$. Then $2^{id_X}=id_{(\MM,\OO)}$;

\smallskip

\noindent(b) Let $(X,\TT)$, $(X',\TT')$, $f$, $\MM$, $\MM'$, $\PP$, $\PP'$, $\OO$ and $\OO'$ be as in Proposition \ref{2^f}, $(X'',\TT'')$ be a topological space, $\MM''=CL(X'')$, $\PP''$ be a subbase for $X''$, $\OO''=\OO_{\PP''}^{\MM''}$, $g:X'\longrightarrow X''$ be a map and $g^{-1}(\PP'')\subseteq\PP'$. Then $2^{g\circ f}=2^g\circ 2^f$.
\end{pro}

\doc (a) Obvious.

(b) Let $f:X\longrightarrow Y$, $g:Y\longrightarrow Z$. Then $2^f(C)=\overline{f(C)}^Y$, $2^g(D)=\overline{g(D)}^Z$ and $2^{g\circ f}(A)=\overline{(g\circ f)(A)}^Z$.
We have that $$(2^g\circ 2^f)(C)=2^g(\overline{f(C)})=\overline{g(\overline{f(C)})}$$ and $$2^{g\circ f}(C)=\overline{(g\circ f)(C)}.$$

So, we have to prove that $\overline{g(\overline{f(C)})}=\overline{(g\circ f)(C)}=\overline{g(f(C))}$, i.e. $\overline{g(\overline{A})}=\overline{g(A)}$, for $A\subseteq X'$. Since $g$ is continuous, this is fulfilled. Thus, $2^{g\circ f}=2^g\circ 2^f$.  \sqs

\begin{pro}\label{vr2^f}
Let $X,Y$ be sets, $\PP\sbe\PP\ap(X)$, $\SS\sbe\PP\ap(Y)$, $\bigcup\PP=X$, $\bigcup\SS=Y$, $\TT$ (resp., $\TT\ap$) be the topology on $X$ (resp., $Y$) generated by the subbase $\PP$ (resp., $\SS$),  $f:X\longrightarrow Y$, $f^{-1}(\SS)\subseteq\PP$, $\MM\subseteq \PP'(X)$, $\MM'=CL(Y,\TT\ap)$,  $\OO=\OO_{\PP}^\MM$ and $\OO'=\OO_{\SS}^{\MM'}$. Let $(X,\TT)$ be a $T_2$-space, $(Y,\TT')$ be a $T_1$-space, $\MM$ be a natural family and $2^f:(\MM,\OO)\longrightarrow(\MM',\OO')$, where $2^f(F)=\overline{f(F)}^Y$, for every $F\in \MM$, be a closed map. Then the map $f$ is closed.
\end{pro}

\doc Let $F\in CL(X)$. Then $\Phi_X(F)=\{\{x\}\st x\in F\}$ is a closed subset of $\Phi_X(X)$ and $\Phi_X(X)$ is a closed subset of $(\MM,\OO)$, i.e. $\Phi_X(F)$ is a closed subset of $(\MM,\OO)$ (see Proposition \ref{vl}). Then $2^f(\Phi_X(F))$ is a closed subset of $(\MM',\OO')$. We have that $2^f(\Phi_X(F))=\{2^f(\{x\})\st x\in F\}=\{\{f(x)\}\st x\in F\}=\Phi_Y(f(F))$. Hence $\Phi_Y(f(F))$ is a closed subset of $\Phi_Y(Y)$. Since $\Phi_Y$ is a homeomorphic embedding, we get that $f(F)$ is a closed subset of $Y$. \sqs

With our next result we generalize a theorem of H.-J. Schmidt \cite[Theorem 11(1)]{Sch}  about {\it commutability between hyperspaces and subspaces} (see also \cite{D1} for similar results).

\begin{theorem}\label{i_A}
Let $(X,\TT)$ be a space, $\PP$ be a subbase for $\TT$, $X\in\PP$, $\OO=\OO_{\PP}^{CL(X)}$. For any $A\subseteq X$, set $\PP_A=\{U\cap A\st U\in\PP\}$ and $\OO_A=\OO_{\PP_A}^{CL(A)}$. Then $i_{A,X}:(CL(A),\OO_A)\longrightarrow(CL(X),\OO)$, where $i_{A,X}(F)=\overline{F}^X$, is a homeomorphic embedding.
\end{theorem}

\doc Let $F$ be a closed subset of $A$, $U\in\PP$ and $U\cap\overline{F}^X\not=\emptyset$ (i.e. $\overline{F}^X\in U^-_{CL(X)}$). Then $U\cap F\not=\emptyset$, i.e. $F\in(U\cap A)^-_{CL(A)}$ and $U\cap A\in\PP_A$. If  $G\in(U\cap A)^-_{CL(A)}$ then $G\cap U\not=\emptyset$. Hence $U\cap\overline{G}^X\not=\emptyset$, i.e. $\overline{G}^X\in U^-_{CL(X)}$. Thus $i_{A,X}((U\cap A)^-_{CL(A)})\subseteq U^-_{CL(X)}$, i.e., using Fact \ref{neprpr},  we get that $i_{A,X}$ is continuous.

Conversely, let $F\in CL(A)$, $U'\in\PP_A$ and $F\cap U'\not=\emptyset$. Then there exists an $U\in\PP$ with $U'=U\cap A$. Hence $U\cap F\not=\emptyset$ and thus $U\cap\overline{F}^X\not=\emptyset$, i.e. $\overline{F}^X\in U^-_{CL(X)}$. If $G\in CL(A)$ and $\overline{G}^X\cap U\not=\emptyset$, then $G\cap U\not=\emptyset$ and hence $G\cap U\cap A\not=\emptyset$, i.e. $G\cap U'\not=\emptyset$. Thus $G\in(U')^-_{CL(A)}$. Hence $i_{A,X}^{-1}(U^-_{CL(X)}\cap i_{A,X}(CL(A)))\subseteq(U')^-_{CL(A)}$. So, using Fact \ref{neprpr},  we get that $i_{A,X}$ is inversely continuous. \sqs

\section{Relations between some topological properties of the spaces and the hyperspaces endowed with a lower-Vietoris-type topology}

In this section we generalize Propositions 3.1 and 3.3 from \cite{CMR}.

\begin{pro}\label{compact}
Let $X$ be a set, $\MM\sbe\PP\ap(X)$, $\MM$ be a natural family, $\OO$ be a lower-Vietoris-type topology on $\MM$ and $(X,\TT_{\OO})$ be a $T_2$-space. Then $(X,\TT_{\OO})$ is compact if and only if $(\MM,\OO)$ is compact.
\end{pro}

\doc Let $(\MM,\OO)$ be compact. From \ref{vl} we have that $\Phi(X)$ is a closed subset of $(\MM,\OO)$. Hence $\Phi(X)$ is compact. Then $(X,\TT_{\OO})$ is compact.

Conversely, let $(X,\TT_{\OO})$ be  compact. Let $\{(U_{\alpha})^-_\MM\st U_{\alpha}\in\PP_\OO,\alpha\in\AA\}$ be a cover of $\MM$. Then $\displaystyle\MM=\bigcup_{\alpha\in\AA}(U_{\alpha})^-_\MM=(\bigcup_{\alpha\in\AA}U_{\alpha})^-_\MM$. We will prove that $\displaystyle\bigcup_{\alpha\in\AA}U_{\alpha}=X$. Indeed, since $\MM$ is a natural family, we get that for every $x\in X$ there exists an $\alpha\in\AA$, such that $\{x\}\in (U_{\alpha})^-_\MM$. Then $x\in U_{\alpha}$. Hence $\displaystyle X=\bigcup_{\alpha\in\AA}U_{\alpha}$. Since $X$ is compact, we get that there exists a finite subcover $\{U_{\alpha_i}\st i=1,\dots,n\}$ of the cover $\{U_{\alpha}\st\alpha\in\AA\}$ of $X$. Then $\displaystyle\bigcup_{i=1}^n(U_{\alpha_i})^-_\MM=(\bigcup_{i=1}^nU_{\alpha_i})^-_\MM=X^-_\MM=\MM$. \sqs

\begin{pro}\label{3.3}
Let $(X,\TT)$ be a topological space, $\MM\subseteq\PP\ap(X)$ be a natural family and $\OO=\OO_\TT^\MM$. Then $w(X,\TT)\leq\tau(\geq\aleph_0)$ if and only if $w(\MM,\OO)\leq\tau$.
\end{pro}

\doc Let $w(\MM,\OO)\leq\tau$. Since $\Phi:(X,\TT)\rightarrow(\MM,\OO)$ is a homeomorphic embedding (see Proposition \ref{vl}), we get that $w(X)\leq\tau$.

Conversely, let $w(X)=\lambda\leq\tau$. Then there exists a base $\BB$ for $X$ with $|\BB|=\lambda$. Then, by Corollary \ref{sl2},
 $\BB^-_{\MM}$ is a subbase for $(\MM,\OO)$. Hence $w(\MM,\OO)\leq|(\BB^-_\MM)^\cap|=|\BB^-_\MM|\leq|\BB|=\lambda\leq\tau$. \sqs

\section{Homotopy, extensions of maps and fixed-point properties in hyperspaces endowed with a low\-er-Vietoris-type topology}

\begin{pro}\label{4.1}
 Let $X$ be a set, $\MM\subseteq\PP\ap(X)$ be a natural family, $\OO$ be a lower-Vietoris-type topology on $\MM$ and
 $X\in\MM$. Then $(\MM,\OO)$ has a trivial homotopy type.
\end{pro}

\doc It is analogous to the proof of Proposition 4.1 from \cite{CMR}. \sqs

\begin{cor}\label{4.2}
Let $X$ be a set, $\MM\subseteq\PP\ap(X)$ be a natural family, $\OO$ be a lower-Vietoris-type topology on $\MM$ and
 $X\in\MM$.
 Then $(\MM,\OO)$ is contractible and locally contractible.
\end{cor}

\begin{pro}\label{4.3}
Let $X$ be a set, $\MM\subseteq\PP\ap(X)$ be a natural family, $\OO$ be a lower-Vietoris-type topology on $\MM$ and
$X\in\MM$. Then $(\MM,\OO)$ is an absolute extensor for the class of all topological spaces (i.e.,  every continuous function $f$ from a closed subspace of a space $Z$ to $(\MM,\OO)$ can be continuously extended to $Z$).
\end{pro}

\doc It is analogous to the proof of  Proposition 4.3 from \cite{CMR}. \sqs

\begin{pro}\label{zatwobv}
Let $(X,\TT)$ be a topological space, $\PP$ be a subbase for $\TT$, $\emptyset\not=\MM\sbe\{X\stm U\st U\in\PP\}$, $\ems\nin\MM$ and $\OO=\OO_\PP^\MM$. Then,
for every $M\in\MM$, $\overline{M}^{\OO}=\{M'\in\MM\st M'\subseteq M\}(=M^+_\MM)$.
\end{pro}

\doc Let $M\in\MM$. If $M=X$ then $\overline{M}^{\OO}=M_{\MM}^+$. So that, let $M\not=X$. It is obvious that $M^+_\MM\subseteq\overline{M}^{\OO}$. Suppose that there exists an $M'\in\MM$ such that $M'\in\overline{M}^{\OO}$ and $M'\setminus M\not=\emptyset$. Let $U=X\setminus M$. Then $U\in\PP$ (because $M\in\MM$) and $M'\cap U\not=\emptyset$. Hence $M'\in U^-_\MM$, but, obviously, $M\not\in U^-_\MM$. Thus $M'\not\in\overline{M}^{\OO}$,  a contradiction. Therefore, $M'\subseteq M$, i.e. $\overline{M}^{\OO}\subseteq M^+_\MM$. \sqs

\begin{theorem}\label{fixed}
Let $(X,\TT)$ be a topological space, $\PP$ be a subbase for $\TT$,  $\ems\neq\MM\sbe\{X\stm U\st U\in\PP\}$, $\ems\nin\MM$ and $X\in\MM$.
If $\OO=\OO_{\PP}^{\MM}$ and $\MM$ is closed under intersections of decreasing subfamilies, then $(\MM,\OO)$ has the fixed-point property.
\end{theorem}

\doc Let $\Psi:(\MM,\OO)\longrightarrow(\MM,\OO)$ be continuous. Set $\KK=\{K\in\MM\st\Psi(K)\subseteq K\}$.

In $\KK$, set $A<B$ if and only if $B\subsetneqq A$. It is obvious that $(\KK,\leq)$ is an ordered set. Since $X\in\KK$ we get that $\KK\not=\emptyset$. We will prove that if $K\in\KK$, then $\Psi(K)\in\KK$. Indeed, let $K\in\KK$. Then $\Psi(K)\subseteq K$.
We have that $\Psi(K)\in\MM$. Then  Proposition \ref{zatwobv} implies that $\Psi(K)\in\overline{K}^{\OO}$. By the continuity of $\Psi$, we get that $\Psi(\overline{K}^{\OO})\subseteq\overline{\Psi(K)}^{\OO}$. Since $\Psi(K)\in\overline{K}^{\OO}$, we obtain that $\Psi(\Psi(K))\in\overline{\Psi(K)}^{\OO}$. Then Proposition \ref{zatwobv} implies that $\Psi(\Psi(K))\subseteq\Psi(K)$. Hence $\Psi(K)\in\KK$.

Let now, $\{K_{\alpha}\st\alpha\in A\}$ be a chain in $(\KK,\leq)$. We have that $\displaystyle\emptyset\not=K=\bigcap_{\alpha\in A}K_{\alpha}\in\MM$. We will prove that $K\in\KK$. Indeed, for every $\alpha\in A$, we have that $K\subseteq K_{\alpha}$. Hence $K\in\overline{K_{\alpha}}^{\OO}$ for every $\alpha\in A$. Since $\Psi(\overline{K_{\alpha}}^{\OO})\subseteq\overline{\Psi(K_{\alpha})}^{\OO}$, we get that $\Psi(K)\in\overline{\Psi(K_{\alpha})}^{\OO}$ for every $\a\in A$. Then, by Proposition \ref{zatwobv},  $\Psi(K)\subseteq\Psi(K_{\alpha})\subseteq K_{\alpha}$,  for every $\alpha\in A$. Hence $\displaystyle\Psi(K)\subseteq\bigcap_{\alpha\in A}K_{\alpha}=K$, i.e. $K\in\KK$.

So, $(\KK,\leq)$ satisfies the hypothesis of the Zorn Lemma. Hence it has a maximal element $K_0$. Since $K_0\in\KK$, we get that $\Psi(K_0)\subseteq K_0$. Suppose that $\Psi(K_0)\not=K_0$. Then $\Psi(K_0)>K_0$, a contradiction.  Therefore, $\Psi(K_0)=K_0$. \sqs

\begin{cor}\label{sl2fix}({\rm \cite{CMR}})
Let $(X,\TT)$ be a compact Hausdorff space.
Then $(CL(X),\UP^-_X)$ has the fixed-point property.
\end{cor}

\begin{cor}\label{sl1fixsp}({\rm \cite{CMR}})
Let $(X,\TT)$ be a compact Hausdorff space and $f:X\lra X$ be a continuous map. Then there exists a compact subspace $K$ of $X$ such that $f(K)=K$.
\end{cor}

\begin{cor}\label{sl3fix}
Let $(X,\TT)$ be a compact Hausdorff space, $\MM\subseteq CL(X)$ be closed under any intersections of decreasing subfamilies and $X\in\MM$. Then $(\MM,\OO_\TT^\MM)$ has the fixed-point property.
\end{cor}

\doc Set $\PP=\TT$ in Theorem \ref{fixed}. \sqs

\begin{cor}\label{sl4fix}
Let $(X,\TT)$ be a continuum, $$\MM=\{K\subseteq X\st K\mbox{ is a non-empty continuum}\},$$ $\PP\supseteq\{X\setminus K\st K\in\MM\}$, $\PP$ be a subbase for $\TT$ and $\OO=\OO_{\PP}^{\MM}$. Then $(\MM,\OO)$ has the fixed-point property.
\end{cor}

\doc By \cite[Corollary 6.1.19]{E2}, $\MM$ is closed under intersections of decreasing subfamilies.  Since $X\in\MM$ and $\emptyset\not\in\MM$, we have  that $(\MM,\OO)$ satisfies the hypothesis of Theorem \ref{fixed}. \sqs

\begin{cor}\label{sl5fix}
Let $(X,\TT)$ be a continuum and $f:X\longrightarrow X$ be a continuous map. Then there exists a continuum $K\subseteq X$, such  that $f(K)=K$.
\end{cor}

\doc Since $f$ is a continuous map, we get, by Proposition \ref{2^f}, that $$2^f:(\MM,\OO_\TT^\MM)\longrightarrow(\MM,\OO_\TT^\MM),$$ where $2^f(K)=f(K)$ and $\MM=\{K\subseteq X\st K\mbox{ is a non-empty continuum }\}$, is a continuous map. Then, by Corollary \ref{sl4fix}, we obtain that there exists a $K\in\MM$ such that $2^f(K)=K$. Hence $f(K)=K$. \sqs

\baselineskip = 0.75\normalbaselineskip

\end{document}